\documentclass[12pt]{amsart}
\usepackage{amssymb} 
\usepackage[mathscr]{eucal}
\usepackage{epsf}

\numberwithin{equation}{section} 

\makeatletter
\def\LaTeX{\leavevmode L\raise.42ex
    \hbox{\kern-.3em\size{\sf@size}{0pt}\selectfont A}\kern-.15em\TeX}
\makeatother

\newcommand{\BibTeX}{{\rm B\kern-.05em{\sc i\kern-.025emb}\kern-.08em\TeX}}

\newtheorem{thm}{Theorem}

\newtheorem{cor}[thm]{Corollary}

\theoremstyle{definition}

\theoremstyle{remark}
\newtheorem{rem}[thm]{Remark} 
\theoremstyle{remark}

\makeatletter
\def\@currentlabel{2.1}\label{e:dispaa}
\def\@currentlabel{2.21}\label{e:dispau}
\def\@currentlabel{2.22}\label{e:dispav}
\def\@currentlabel{2.23}\label{e:dispaw}
\def\@currentlabel{2.24}\label{e:dispax}
\def\theequation{\thesection.\@arabic\c@equation}
\makeatother

\makeatletter
\def\alphenumi{%
  \def\theenumi{\alph{enumi}}%
  \def\p@enumi{\theenumi}%
  \def\labelenumi{(\@alph\c@enumi)}}
\makeatother

\newcommand{\marginnote}[1]
{
}

\newcounter{gm}

\newcounter{bk}

\begin{document}
\hfill{To appear in {\it  J. of Math. Fluid Mech.}}

\bigskip
\bigskip

\centerline{\it To Victor Yudovich on the occasion of his 70th birthday}
\bigskip
\bigskip
\bigskip

\title[Asymptotic directions, Monge-Amp\`ere and diffeomorphism groups]
{Asymptotic directions, Monge-Amp\`ere equations and the geometry 
of diffeomorphism groups} 
\author{Boris Khesin}
\address{B.K.: Department of Mathematics, University of Toronto, 
ON M5S 3G3, Canada} 
\thanks{The first author was partially supported by an NSERC grant.}
\email{khesin@math.toronto.edu} 
\author{Gerard Misio\l ek}
\address{G.M.: Department of Mathematics, University of Notre Dame, 
IN 46556, USA}
\email{gmisiole@nd.edu} 

\begin{abstract}
In this note we obtain the characterization for asymptotic directions on 
various subgroups of the diffeomorphism group. We give a simple proof 
of non-existence of such directions for area-preserving diffeomorphisms 
of closed surfaces of non-zero curvature.
Finally, we exhibit the common origin of the Monge-Amp\`ere equations
in 2D fluid dynamics and mass transport.
\end{abstract}

\maketitle 


\section{Asymptotic directions on diffeomorphism groups}
\smallskip

The group of volume-preserving diffeomorphisms of a Riemannian 
manifold plays a fundamental role in the geometrical study of
the Euler equation of hydrodynamics on the manifold \cite{ar}. 
In this paper we consider another equation, 
the Monge-Amp\`ere equation, and discuss its universality 
in the context of diffeomorphism groups. This equation occurs in two main 
contexts: as the equation of asymptotic directions in 2D hydrodynamics
and in the optimal transport problem in any dimension.

\subsection{Characterization of asymptotic directions on subgroups} 
Let $M$ be a compact $n$-dimensional manifold without boundary and equipped 
with a Riemannian metric $g$. Let $\mathcal{D}(M)$ be (the connected 
component of the identity in) the group of all diffeomorphisms 
of $M$. Its  tangent space at the identity diffeomorphism consists of 
smooth vector fields on $M$.  
The tangent space at a point $\eta$ consists of vector fields 
``reparameterized by $\eta$,'' i.e. of the maps 
$X_\eta : M \to TM$ with $X_\eta(x) \in T_{\eta(x)}M$. 
Define the "flat" $L^2$-metric on $\mathcal{D}(M)$ by assigning to 
each tangent space the inner product 
\begin{equation} \label{L2} 
g_\eta(X_\eta, Y_\eta) 
:= 
\int_M g_{\eta(x)}\big( X_\eta(x),Y_\eta(x) \big) \; dV(x),
\end{equation} 
where $dV$ denotes the Riemannian volume form on $M$. 

Let $\mathcal{SD}(M)$ be the subgroup of volume-preserving diffeomorphisms.
The restriction of the $L^2$-metric to this subgroup is right-invariant
and of particular importance
in hydrodynamics. In  \cite{ar} Arnold
showed that geodesics in $\mathcal{SD}(M)$ correspond to 
motions of an ideal fluid in $M$ described by the Euler equations 
$$
\partial_t X + \nabla_X X = - \nabla p, 
\quad 
\mathrm{div}\, X = 0 
$$
for the fluid velocity field $X$. 

As shown by Ebin and Marsden in \cite{em}
the group $\mathcal{D}^s(M)$ of all diffeomorphisms of Sobolev class $H^s$ 
(as well as its various subgroups) 
can be viewed for $s>n/2 + 1$ as an infinite dimensional Hilbert manifold. 
All of the arguments in this paper can be rigorously developed in 
the Sobolev framework. However, to present the geometric ideas 
we will keep things formal and drop the index $s$ in what follows.

While the Euler equation depends only on the intrinsic Riemannian geometry 
of $\mathcal{SD}(M)$, it is also of interest to study its exterior geometry
as a Riemannian submanifold in $\mathcal{D}(M)$. 
In particular, one can consider asymptotic directions in $\mathcal{SD}(M)$. 
A vector tangent to a Riemannian submanifold is
{\it asymptotic} if the geodesics issued in the direction of this vector, one 
in the submanifold and the other in the ambient manifold, 
have a second order of tangency. (Note that in general two geodesics 
with a common tangent will have only a simple, 
i.e. first order, tangency.) 
More formally, asymptotic vectors are singled out  by the condition that 
the second fundamental form evaluated on these vectors is zero. 
A curve whose tangent is asymptotic at each point 
is called an asymptotic line. An asymptotic line 
is a geodesic in the submanifold if and only if 
it is also a geodesic in the ambient manifold.

A  description of asymptotic vectors in the group of 
volume-preserving diffeomorphisms 
was  given by Bao and Ratiu in \cite{br},  
while asymptotic geodesics in $\mathcal{SD}(M)$ were  studied 
in \cite{m1} (as pressure-constant flows). 
We begin with the following convenient characterization of 
these vectors. 

\begin{thm} \cite{br}\label{thm:br}
A  vector field $X$ on a manifold $M$ is an  
asymptotic direction  for $\mathcal{SD}(M)$ at the identity 
diffeomorphism if and only if 
\begin{equation}\label{div}
\mathrm{div}\,\nabla_{X}X = \mathrm{div}\,X = 0. 
\end{equation}
\end{thm}

Similarly, a vector $X_\eta$ is asymptotic to  $\mathcal{SD}(M)$ at the 
diffeomorphism $\eta$ if and only if its right translation
$X:=X_\eta\circ \eta^{-1}$ is asymptotic to $\mathcal{SD}(M)$ at the identity.
If $M$ is two-dimensional then more can be said: 

\begin{cor}\cite{br}
The stream function $\psi$ of an asymptotic vector 
field $X$ in 2D satisfies a Monge-Amp\`ere 
equation 
\begin{equation} \label{MA} 
\det [D^2\psi]=\frac{g\cdot K}{2} |\nabla\psi|^2, 
\end{equation}
where $g=\det(g_{ij}),~~ K$ is the Gaussian curvature function of $M$ 
and $\det [D^2\psi]$ is the Hessian of $\psi$ for the field $X=J\nabla\psi$
with respect to the symplectic area form $\omega=dV$ on $M$.
\end{cor}

Below we give a characterization of asymptotic 
vectors in a general setting. Let $\mathcal B$ be a subgroup
of $\mathcal{D}(M)$ and let $\mathfrak b$ denote its Lie algebra. 
We assume that $\mathfrak b$ is a closed subspace, and therefore it 
has an orthogonal 
complement with respect to the $L^2$ inner product (\ref{L2}). 

\begin{thm} \label{gen}
A vector field $X$ on $M$ is asymptotic for $\mathcal B$ 
if and only if 
\begin{equation}\label{eq:gen}
X \in \mathfrak b 
\quad 
\mathrm{and} 
\quad 
\nabla_{X}X \in \mathfrak b\,.
\end{equation}
\end{thm}

If $\mathcal{B}=\mathcal{SD}(M)$ then the Lie algebra 
$\mathfrak b$ consists of divergence-free vector fields 
and we recover Theorem \ref{thm:br}. 
Moreover, if $M$ has a boundary  then 
the diffeomorphisms from the subgroup $\mathcal B$ leave it 
invariant and hence the fields from the subalgebra $\mathfrak b$ 
are tangent to the boundary of $M$. 

\begin{cor}\cite{br} If $M$ has a boundary $\partial M$, 
then $X$ is asymptotic 
for $\mathcal{SD}(M)$ if in addition to equations (\ref{div}) 
$X$ satisfies the conditions 
\begin{equation} \label{BC} 
g(X,n)=g(\nabla_{X}X,n)=0
\end{equation} 
where $n$ is the normal to $\partial M$. 
\end{cor} 


Suppose next that $M$ is a symplectic manifold of dimension $2n$ equipped with 
a symplectic 2-form $\omega$. 
\begin{cor} {\bf a)}
Let $\mathcal{B}=Symp(M)$ be the subgroup of symplectic diffeomorphisms 
with the corresponding Lie subalgebra $\mathfrak{b}=symp(M)$ of vector 
fields preserving $\omega$. 
Then a vector field $X$ is asymptotic if and only if 
$$
L_X \omega = 0 
\quad 
\mathrm{and} 
\quad 
L_{\nabla_{X}X}\omega =0.
$$ 

{\bf b)} Let $\mathcal{B} = Ham(M)$ be the (generally speaking, smaller) 
subgroup of Hamiltonian diffeomorphisms. Then a vector field $X$ on $M$ 
is asymptotic if and only if 
$$ 
\omega(X,\cdot) = d\psi 
\quad 
\mathrm{and} 
\quad 
\omega(\nabla_{X}X, \cdot ) = d\phi 
$$ 
for some smooth functions $\psi$ and $\phi$ on $M$. 
\end{cor} 
If $\dim M=2$ any divergence-free field is locally Hamiltonian and 
the latter equation is rewritten in (\ref{MA}) 
as the Monge-Amp\`ere equation on its Hamiltonian, or stream function. 

Curiously, this Monge-Amp\`ere property does not survive when passing to
Hamiltonian fields  in higher dimensions, unlike what was conjectured at the end
of \cite{ak2}. Indeed, already in the flat 4-dimensional case one obtains a 
system of three equations on pairwise products of second derivatives, while
the corresponding Monge-Amp\`ere with $\det [D^2\psi]$ would include the four-term 
products of second derivatives.

Another interesting example is provided by a contact manifold $M$ with 
a contact structure
(i.e. maximally non-integrable distribution of hyperplanes) $\tau$.

\begin{cor}
Let $(M,\tau)$ be a contact manifold, and $\mathcal B$ a subgroup
consisting of contact diffeomorphisms. Then the vector field $X$ is asymptotic to
$\mathcal B$ if and only if 
$$
both~~X~~and~~\nabla_X X ~~are~~contact.
$$
\end{cor}

If $\alpha$ denotes a 1-form defining the contact structure $\tau$ then 
the Lie algebra $\mathfrak b:= \{X: L_X\alpha = f\alpha\}$. 
For a contact vector field
$X$ consider its contact Hamiltonian function $K_X:=\alpha(X)$ with respect to the
form $\alpha$. By  rewriting the condition 
$\nabla_X X \in \mathfrak b$ for asymptotic vectors on $\mathcal B$ one can obtain
an analog of the Monge-Amp\`ere equation (\ref{MA})
on the Hamiltonian in the contact 3-dimensional case.

\medskip

\subsection*{Proof of Theorem \ref{gen}} 
The Lie algebra of $\mathcal{D}(M)$ splits orthogonally into 
$\mathfrak{b} \oplus \mathfrak{b}^\perp$ 
inducing a corresponding splitting at any point in $\mathcal B$ 
by right translations. 
Let $\bar\nabla$ denote the smooth Levi-Civita connection of the $L^2$-metric 
(\ref{L2}) on $\mathcal{D}(M)$ (see Ebin and Marsden \cite{em}). Let $X$ 
and $Y$ be two elements of $\mathfrak b$ and extend them to 
right-invariant vector fields $X_\eta = X\circ\eta$ and 
$Y_\eta = Y\circ\eta$ on $\mathcal B$. 
Decomposing into unique tangential and normal components 
(Gauss equation) we obtain 
$$
\bar\nabla_{X_\eta} Y_\eta 
= 
P_\eta(\bar\nabla_{X_\eta}Y_\eta) + (\bar\nabla_{X_\eta}Y_\eta)^\perp. 
$$
where $P_{\eta}$ is the projection onto the tangent space to $\mathcal B$ 
at the point $\eta$. 
The tangential component can be used to define a smooth right-invariant 
connection on $\mathcal B$ while the normal component 
$(\bar\nabla_{X_\eta}Y_\eta)^\perp$ is interpreted 
as the second fundamental form of $\mathcal B$ in $\mathcal{D}(M)$. 
Since 
$$
\bar\nabla_{X_\eta}Y_\eta = (\nabla_XY)\circ\eta 
$$
we see that at the identity the $L^2$ covariant derivative $\bar\nabla$ 
is given by $\nabla_XY$. 
Recall that a vector $X$ in $\mathfrak b$ is asymptotic if 
the second fundamental form evaluated at $X$ is zero. 
This implies that $(\nabla_XX)^\perp = 0$ and Theorem \ref{gen} follows. 
$\qquad\qquad\qquad\qquad\quad\square$

\begin{rem} 
One can see how the above general consideration works in the case 
$\mathcal{B}=\mathcal{SD}(M)$. 
First, recall that an arbitrary vector field $X$ on $M$ can be decomposed 
into $L^2$-orthogonal divergence-free and gradient parts 
$$
X = P_{id}( X) + X^{\nabla}
$$
where $P_{id}$ is now the projection onto the divergence-free fields, and 
$X^{\nabla}$ is the gradient part of $X$: explicitly
$X^{\nabla}= \nabla\Delta^{-1}\mathrm{div}\, X$.
By right invariance this induces a corresponding splitting at any point 
in $\mathcal{SD}(M)$. The tangential component of $\bar\nabla$ defines a smooth 
right-invariant connection on $\mathcal{SD}(M)$. Its normal component 
is the second fundamental form of the subgroup 
$\mathcal{SD}(M)$ in $\mathcal{D}(M)$.
Thus, for any right-invariant vector fields $X_\eta = X\circ\eta$ and 
$Y_\eta=Y\circ\eta$ on $\mathcal{SD}(M)$ we again have the Gauss equation
\begin{equation} \label{gauss} 
\bar\nabla_{X_\eta}Y_\eta 
= 
P_\eta (\bar\nabla_{X_\eta}Y_\eta) + s(X_\eta, Y_\eta) 
\end{equation} 
where 
$$
\bar\nabla_{X_\eta}Y_\eta 
= 
( \nabla_X Y )\circ\eta 
\quad 
\mathrm{and} 
\quad 
s(X_\eta ,Y_\eta ) 
= 
\big(\nabla_X Y\big)^\nabla\circ\eta. 
$$
Since a vector $X$ tangent at the identity $\eta=id$ is asymptotic if $s(X,X) = 0$ 
we immediately obtain the characterization in Theorem \ref{thm:br}. 
\end{rem} 

\begin{rem} Asymptotic geodesics (or pressure-constant flows) are 
of interest in the Lagrangian approach to hydrodynamic stability theory. 
Typical examples are plane-parallel flows on the flat torus $\mathbb{T}^2$. 
These flows can be considered unstable in the following sense. 
For any such flow it can be shown that sectional curvatures of 
$\mathcal{SD}(\mathbb{T}^2)$ along the corresponding geodesic are always 
non-positive. Therefore, by a suitable variant of the Rauch comparison 
argument, all linear perturbations in Lagrangian coordinates (Jacobi fields 
along the geodesic) must grow at least linearly in time. 
On the other hand, there are flows with non-positive curvatures for which 
one can show that the growth must be at most polynomial 
(see Preston \cite{p} for details). 
\end{rem} 




\subsection{Non-existence of asymptotic directions}
Asymptotic directions are not always in good supply. 
In \cite{br} it is shown that if $M$ is a two-dimensional compact 
surface of revolution without boundary then any axially symmetric smooth 
solution of the Monge-Amp\`ere equation is constant away from the 
cylindrical (i.e. fixed radius) bands of the surface.
The strongest result in this direction so far 
states that for a compact closed surface $M$ of 
positive curvature there are no asymptotic directions 
on $\mathcal{SD}(M)$, see \cite{pa}. 
A similar result holds for surfaces with boundary. (Palmer also showed 
that there is no direct analog of this result in higher dimensions: 
e.g. every left-invariant vector field on a compact Lie group equipped 
with a bi-invariant metric solves the equations (\ref{div}), 
i.e. is asymptotic. 
This had been observed previously for the three sphere $S^3$ in \cite{m1}.)

Here we prove the following generalization of Palmer's result. 

\begin{thm}
If $M$ is a compact closed surface of nowhere zero curvature $K$, 
then the Monge-Amp\`ere equation (\ref{MA}) 
admits no non-constant solutions. In particular, in this case 
$\mathcal{SD}(M)$ has no asymptotic directions. 
\end{thm} 
Note that the surface can be of any non-zero Euler characteristic 
(the case of the torus 
is ruled out by the Gauss--Bonnet theorem).
\begin{proof} 
Recall that for any vector field $X$ on $M$ we have the identity 
\begin{equation} \label{nabla} 
\mathrm{div}\nabla_X X 
= 
r(X,X) + \mathrm{tr}(D X)^2 
+ 
L_X (\mathrm{div}\, X) 
\end{equation} 
where $r$ denotes the Ricci curvature of the metric $g$ and $L_X$ is the 
Lie derivative along $X$ (see for example \cite{tay}). 
If $X$ is divergence-free then the last term on the right side of 
(\ref{nabla}) drops out. 

Consider the length function $f:= g(X,X)$. Since $M$ is compact $f$ must 
attain a maximum at some point $x_0$. 
Choosing normal coordinates at that point with 
$g_{ij}(x_0)=\delta_{ij}$ and $\partial g_{ij}/\partial x^k(x_0)=0$, 
we obtain 
$$
0=df(x_0) 
= 
2\sum_{jk} X^j(x_0)\, \frac{\partial X^j}{\partial x^k}(x_0)\, dx^k \,.
$$
Since at the point where $f$ has a maximum we must also have $X(x_0)\neq 0$, 
this implies that the Jacobi matrix $DX$ is degenerate at $x_0$. 
Therefore, rearranging the terms of the trace and using the fact 
that the divergence of $X$ is zero, we get 
$$
\mathrm{tr}(D X)^2(x_0) 
= 
\sum_{ij} 
\frac{\partial X^j}{\partial x^i}(x_0) \frac{\partial X^i}{\partial x^j}(x_0) 
= 
-2 \det{[DX(x_0)]}=0. 
$$
Substituting into the relation (\ref{nabla}) we find that 
$$
\mathrm{div}\, \nabla_X X(x_0) 
= K(x_0) \, g(X, X)(x_0)\,, 
$$
because in two dimensions the Ricci and Gaussian curvatures coincide.
However, if $X$ is asymptotic this implies that 
$$ 
0 = K(x_o) \, g(X,X)(x_o) 
$$ 
contradicting the assumption $K\neq 0$ on $M$, and in particular 
at the point $x_0$. 
\end{proof} 

\begin{rem}
Note that 
the proof above works for a $C^1$ field $X$, thus improving on 
the $C^2$ assumption used in \cite{pa}.
\end{rem} 
The following 
is a similar result 
for the case of a manifold with boundary:

\begin{thm}
Let $M$ be a compact surface of nowhere zero curvature $K$ 
with smooth boundary. 
Assume that the geodesic curvature $k_g$ of $\,\partial M$ vanishes, 
at most, at finitely many points. 
Then the Monge-Amp\`ere equation (\ref{MA}) along with the boundary 
condition (\ref{BC}) admits no non-constant solutions.
\end{thm} 

\begin{proof}
The beginning of the proof follows \cite{pa}. Assuming that such 
an $X$ exists we show that $X$ has to vanish on the boundary. 
Indeed, for any point $x\in \partial M$ where $X(x)\not=0$ we have that 
$X$ is tangent to $\partial M$ and 
$$ 
0=(\nabla_X X,n)(x)=k_g(x)\cdot g(X,X)(x)\,. 
$$ 
The latter shows that the geodesic curvature $k_g$  must vanish at $x$. 
By the assumption, this implies that $X$ can be non-zero only at finitely 
many points of the boundary. Hence it is identically zero on $\partial M$ 
by continuity.

The proof of the boundary-free case is now applicable: 
the function $f:= g(X,X)$ is also zero on the boundary 
$\partial M$ and must therefore attain a maximum in the interior of $M$.
\end{proof}


\section{Mass transport and foliations on diffeomorphism groups} 

\subsection{Monge-Amp\`ere equation in optimal mass transport}
A somewhat different Monge-Amp\`ere equation arises in the theory
of mass transport. Namely, let $d\mu(x)=m(x)dV(x)$
and $d\nu(y)=n(y)dV(y)$ be two smooth measures on a manifold $M$ 
of any dimension.
The Jacobian of a map $\eta$ which sends the measure $\mu$ to 
the other $\nu=\eta_* \mu$ satisfies the relation
\begin{equation}\label{jac}
n(\eta(x))\det[D\eta(x)]=m(x).
\end{equation}
The property of optimality means that the 1-parameter family of maps 
$\eta_t(x)$ describes a geodesic curve on the space of densities
with respect to the Wasserstein $L^2$-metric.
The latter is the `transport' $L^2$ metric on the space of densities:
the distance between two densities on $M$ is the cost
of transporting one of them to the other with the $L^2$-cost function.


If $M$ is a domain in $\mathbb{ R}^k$
one can see that an optimal map $\eta: M\to M$ 
has to be the gradient $\eta =\nabla \phi$
of a convex function $\phi$, 
i.e. the equation on the potential $\phi$ assumes the
Monge-Amp\`ere form (see, e.g, \cite{b}):
\begin{equation}\label{transp}
\det[D^2\phi]=\frac{m(x)}{n(\nabla\phi(x))}.
\end{equation}
In the case of an arbitrary manifold $M$ the optimality condition is
expressed in terms of convexity of the function $\phi$ with respect 
to the metric on $M$, see \cite{mc}. The corresponding Monge-Amp\`ere 
equation has the same form as above.



\subsection{Universality of the (pre-) Monge-Amp\`ere equation}
Although the two 
Monge-Amp\`ere equations discussed above look rather different 
(equation (\ref{MA}) is for asymptotic {\it vector fields}, while 
equation (\ref{transp}) is for potentials of optimal {\it diffeomorphisms}),
both have a common origin. It turns out that they 
can be viewed as projections of 
the {\it dispersionless Burgers equation}: 
\begin{equation} \label{B} 
\partial_t X + \nabla_X X = 0\,,
\end{equation} 
which can be thought of as a {\it pre}-Monge-Amp\`ere equation.
This equation describes geodesics with respect to 
the "flat" $L^2$-metric (\ref{L2})  on the group of all 
diffeomorphisms $\mathcal{D}(M)$ for $M$ of any dimension.
According to this equation each fluid particle moves along a geodesic in $M$. 

Note that there is a natural fibration on $\mathcal{D}(M)$ given by 
the projection $\pi$ onto densities ${\mathcal P}(M)$. 
Two diffeomorphisms belong to the same fiber if they move a fixed
density (say, the constant density $O$ for a compact $M$) 
to one and the same density. 
In particular, for a compact $M$ one considers 
densities with the same total mass, and
$\mathcal{SD}(M)=\pi^{-1}({O})$ is one of the fibres. 
This projection is a Riemannian submersion onto ${\mathcal P}(M)$ 
equipped with the Wasserstein (see \cite{o}).

Consider the ``horizontal" geodesics in $\mathcal{D}(M)$, which  
project to geodesics in ${\mathcal P}(M)$, see Fig.1. 
This projection means that instead of solutions of the
Burgers equation (\ref{B}) we are interested only in the divergence 
of the corresponding fields, i.e. in how they act on densities.
The latter is described by the
infinitesimal version of the Monge-Amp\`ere equation
(\ref{transp}). The flow of this infinitesimal
version delivers a solution to the Monge-Amp\`ere equation
for the potential of a gradient diffeomorphism, see \cite{b}.


\begin{figure}[htb]
\centerline{\epsfbox{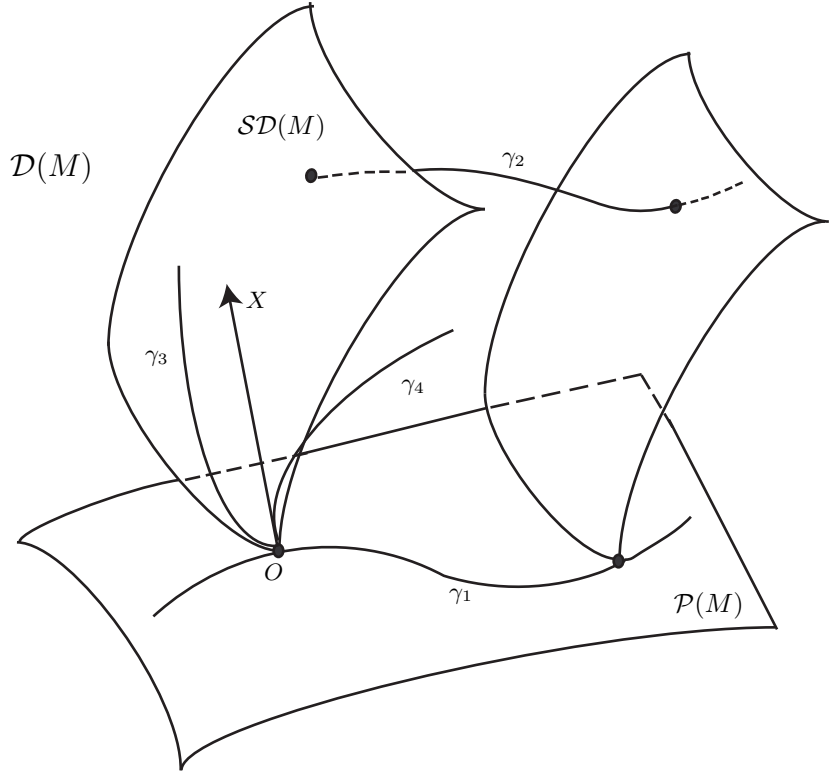}}
\begin{picture} (0, 34)
\put(-70,276){\Small $\mathcal{SD}(M)$} \put(-60,107){\SMALL $O$} 
\put(95,95){\Small $\mathcal{P}(M)$} 
\put(10,100){\SMALL $\gamma_1$} 
\put(30,265){\SMALL $\gamma_2$} 
\put(-105,190){\SMALL $\gamma_3$} 
\put(-8,179){\SMALL $\gamma_4$}
\put(-67,210){\SMALL $X$}  
\put(-160,260){ $\mathcal{D}(M)$}
\end{picture}
\vskip -0.5in
\caption{$\gamma_1$ is a geodesic in the space of densities ${\mathcal P}(M)$,
$\gamma_2$ is a horizontal lift of $\gamma_1$ to ${\mathcal D}(M)$,
$\gamma_3$ is an asymptotic (or `vertical') geodesic in  ${\mathcal SD}(M)$,
$\gamma_4$ is an `almost vertical' geodesic for  
${\mathcal SD}(M)\subset {\mathcal D}(M)$ with the asymptotic direction $X$.}
\label{fig} 
\end{figure}


On the other hand, asymptotic directions on $\mathcal{SD}(M)$
correspond to  ``almost vertical'' geodesics in the space $\mathcal{D}(M)$ with
respect to the same flat $L^2$ metric.
More precisely, an asymptotic {\it geodesic} (given by the Burgers equation)
is ``vertical'' since it joins  diffeomorphisms belonging to the same 
fiber $\mathcal{SD}(M)$, i.e. it joins
different pre-images of the point $O$ on  the quotient.
Asymptotic {\it directions} are  vertical initial vectors which
correspond to geodesics having second order of tangency with the fiber.
In other words, the projections of these geodesics to the base leave the 
initial point 
$O$ very slowly. Indeed, recalling that taking the projection is the same as 
taking the gradient part, we identify here both of equations (\ref{div}): 
along with
the Burgers equation they imply that $\mathrm{div} X=\mathrm{div} X_t=0$.
These relations mean that the projection of the corresponding geodesic 
to the  base ${\mathcal P}(M)$ is a curve which starts at the point $O$ 
with  zero velocity and zero acceleration, respectively. 
(For non-asymptotic directions
only the velocity is zero, while the acceleration is not, since 
typical geodesics have the first order of tangency to  $\mathcal{SD}(M)$.)

The above discussion can be summarized in the following statement:
~
\begin{thm} 
The Monge-Amp\`ere equation (\ref{MA}) in 
2D fluid dynamics and (the infinitesimal form of) the Monge-Amp\`ere 
equation in mass transport (\ref{transp})
are projections of the same dispersionless Burgers equation.
\end{thm} 

\begin{rem}
Yet another appearance of the  Monge-Amp\`ere equation is related to the discrete 
version of the 2D Euler equation for a special energy functional, due to
Moser and Veselov \cite{mv}. A solution of the 
discretized Euler equation for an ideal fluid is a recursive sequence of
diffeomorphisms. The Monge-Amp\`ere equation provides the constraint on the initial
condition ensuring that all diffeomorphisms of the discrete geodesics 
are area-preserving. This restriction is similar to the Monge-Amp\`ere condition
singling out asymptotic directions among 2D divergence-free vector fields.
\end{rem}


\subsection{Relations of diffeomorphism groups}
One of the most interesting problems in the context of diffeomorphism groups 
is the existence of a shortest path or a geodesic between any two diffeomorphisms.
This question is interesting already in dimension $n=1$, for the group
$\mathcal{D}(S^1)$ equipped with either $L^2$ or $H^1$ metric. 
Locally such a geodesic always exists (see \cite{m2}). 

In finite-dimensional geometry such questions are answered using 
the Hopf-Rinow theorem, but in infinite dimensions things become complicated. 
Grossman and Atkin  constructed examples which 
show that even complete infinite-dimensional Hilbert Riemannian 
manifolds have points that cannot be joined by any geodesic. 
For instance, for the groups $\mathcal{SD}(M^n)$ with $n\geq 3$ 
Shnirelman \cite{sh} constructed examples of pairs of diffeomorphisms 
for which a shortest path does not exist.

In the study of this problem the following relation of 
the groups of diffeomorphisms of one- and two-dimensional manifolds can be 
useful. Consider the fibration of  all diffeomorphisms of a surface $M$
over the space of densities on $M$, 
which is exactly the context  of optimal mass transport:
$$ 
\mathcal{SD}(M) 
\hookrightarrow 
\mathcal{D}(M) 
\xrightarrow{~} 
C^\infty(M)\,.
$$ 
Here the left arrow is an isometric embedding with respect to the corresponding 
$L^2$-metrics discussed above, and the right arrow  is
a Riemannian submersion into the space of densities $C^\infty(M)$
equipped with the Wasserstein $L^2$-metric. 

Now, consider a surface $M^2$ whose boundary is the circle $S^1$. 
Similarly to the above, we have the following fibration 
$$ 
\mathcal{SD}_o(M^2) 
\hookrightarrow 
\mathcal{SD}(M^2) 
\xrightarrow{~} 
\mathcal{D}(S^1) 
$$ 
where $\mathcal{SD}_o(M^2)$ denotes the subgroup of smooth area-preserving 
diffeomorphisms of the surface that are pointwise fixed on the boundary $S^1$.


It would be interesting to study the Riemannian properties of this fibration.
In this context the shortest path problem in the group 
$\mathcal{D}(S^1)$ looks like the problem for an optimal path in
the group $\mathcal{SD}(M^2)$ or $\mathcal{D}(M^2)$,
assuming that the projection to $\mathcal{D}(S^1)$ is a Riemannian submersion.
The latter is a problem of optimal transport with a prescribed 
boundary map: we are connecting in an optimal way two disk diffeomorphisms 
which are liftings of the given diffeomorphisms of the circle. 
Here one could hope to  employ the following  reasoning.
While the optimal map exists and is unique between any two convex domains,
it automatically determines how the boundary is mapped. Hence one can 
almost never solve the same problem of finding the optimal map and 
simultaneously satisfy a particular boundary condition. This would allow one 
to conclude that the shortest path does not exist for almost all pairs of 
circle diffeomorphisms.

\bigskip


\end{document}